\documentclass[12pt]{amsart}
\usepackage{amssymb}
\usepackage{amsfonts}
\usepackage{latexsym}
\usepackage{amscd}
\usepackage[mathscr]{euscript}
\usepackage{xy} \xyoption{all}

\newcommand{\Z}{\mathbb{Z}}

\newcommand{\N}{\mathbb{N}}
\newcommand{\ord}{\mathrm{ord}}

\theoremstyle{plain}
\newtheorem{thm}{Theorem}
\newtheorem{lem}[thm]{Lemma}
\newtheorem{prop}[thm]{Proposition}

\theoremstyle{definition}

\newtheorem{example}[thm]{Example}

\newcount\parracount
\newcount\recount
\usepackage{setspace}

\title[Matrix type of purely infinite simple Leavitt path algebras]{The matrix type of purely infinite simple \\ Leavitt path algebras}

\begin{document}

\author{Gene Abrams and Christopher Smith}
\address{
Department of Mathematics, University of Colorado \\
Colorado Springs CO 80918 U.S.A.}
\email{abrams@math.uccs.edu \ \ cdsmith@gmail.com}
\thanks{The first  author is partially supported by the U.S. National Security Agency
 under grant number H89230-09-1-0066.} 

\maketitle

\begin{abstract}
Let $R$ denote the purely infinite simple unital Leavitt path algebra $L(E)$.   We
completely determine the pairs of positive integers $(c,d)$ for which there is
an isomorphism of matrix rings
${\rm M}_c(R)\cong {\rm M}_d(R)$, in terms of the order
of $[1_R]$ in the Grothendieck group $K_0(R)$.
\end{abstract}

\doublespacing

For a row-finite directed graph $E$ and field $k$, the Leavitt path algebra $L_k(E)$ has been
defined in \cite{AA1} and \cite{AMP}, and further investigated in numerous subsequent
articles.   Purely infinite
simple rings were introduced in \cite{AGP};  the purely infinite simple Leavitt path
algebras were explicitly described in \cite{AA2}.  All terminology used in this article can be found
in these four references. We denote $L_k(E)$ simply by  $L(E)$ throughout.

In this short note we present necessary and sufficient conditions for the existence of a ring
isomorphism between the matrix rings ${\rm M}_c(L(E))$ and ${\rm M}_d(L(E))$ (thereby yielding the so-called
{\it Matrix Type} of $L(E)$), whenever $L(E)$
is both purely infinite simple and unital.  ($L(E)$ is unital precisely when the graph $E$ is finite.)
The sufficiency of these conditions utilizes the deep ``algebraic Kirchberg Phillips Theorem"
\cite[Theorem 2.5]{ALPS} for Leavitt path algebras: If $L(E)$ and $L(F)$ are Morita equivalent
purely infinite simple unital rings, and there exists an isomorphism $\varphi: K_0(L(E))
\rightarrow K_0(L(F))$ for which
$\varphi([1_{L(E)}]) = [1_{L(F)}]$, then $L(E) \cong L(F)$.


The following result is well-known, but we prove it here for completeness.

\begin{lem}\label{grouptheoryresult}
Let $G$ be a finitely generated abelian group (written additively).  Let
 $x \in G$ be an element of finite order $n$, and let
$c, d \in \N$.  There exists an automorphism $\varphi : G \to G$ with $\varphi(cx) = dx$
if and only if $\gcd(c,n) = \gcd(d,n)$.
\end{lem}
\begin{proof}
($\Leftarrow$)  \ Since $G$ is finitely generated, $G \cong \Z_{{p_1}^{k_1}} \oplus \cdots \oplus
\Z_{{p_s}^{k_s}} \oplus \Z^t$ for some nonnegative integers $s,t$, (not necessarily distinct)
primes $p_i$ ($1\leq i \leq s$), and $k_i \in \N$.  Since $x$ has finite order, we have
$x = (x_1, \ldots, x_s, 0, \ldots, 0)$ with $x_i \in \Z_{{p_1}^{k_1}}$.   Let
$m_i$ denote $\ord(x_i)$.  Then $m_i|n$, so we have  $\gcd(c,m_i)
= \gcd(\gcd(c,n),m_i)$, which by hypothesis equals $\gcd(\gcd(d,n),m_i)$, which
 (again using $m_i|n$) equals $\gcd(d,m_i)$.  Consequently, $\ord(cx_i) = m_i / \gcd(c,m_i)
= m_i / \gcd(d,m_i) = \ord(dx_i)$.  Since $\Z_{{p_i}^{k_i}}$ is cyclic, and $cx_i$ and $dx_i$
have the same order in $\Z_{{p_i}^{k_i}}$, there exists an automorphism $\varphi_i$ of  $\Z_{{p_i}^{k_i}}$ with $\varphi_i(cx_i) = dx_i$.
Now define $\varphi = \varphi_1 \oplus \cdots \oplus \varphi_s \oplus \mathrm{id}_\Z
\oplus \cdots \oplus \mathrm{id}_\Z \in {\rm Aut}(G)$; then clearly  $\varphi(cx) = dx$, as desired.

($\Rightarrow$)  \ Conversely, suppose $cx \stackrel{\tiny \varphi}{\mapsto} dx$ for some automorphism $\varphi$ of $G$.
Then $n / \gcd(c,n) = \ord(cx) = \ord(dx) = n / \gcd(d,n)$, so $\gcd(c,n) = \gcd(d,n)$.
\end{proof}

Our first of two main results generalizes to all purely infinite simple unital Leavitt path
algebras $L(E)$ a result known previously for the Leavitt algebras $L_q$.   We note for later use that when $E$ is finite, the semigroup $\mathcal{V}^\ast(L(E))$, and therefore the group $K_0(L(E))$, is
finitely generated by \cite[Theorem 3.5]{AMP}.

\begin{thm}\label{regularmodulefiniteorder}
Let $E$ be a graph for which $L(E)$ is purely infinite simple unital. Suppose $[1_{L(E)}] \in K_0(L(E))$
has finite order $n$.  Then, for any $c, d \in \N$, there exists an isomorphism of
matrix rings ${\rm M}_c(L(E)) \cong {\rm M}_d(L(E))$ if and only if $\gcd(c,n) = \gcd(d,n)$.
\end{thm}
\begin{proof}
($\Rightarrow$) \  Because $K_0(L(E))$ is finitely generated, Lemma \ref{grouptheoryresult}
ensures that there exists $\varphi \in \mathrm{Aut}(K_0(L(E)))$ having
$\varphi(c[1_{L(E)}]) = d[1_{L(E)}]$.

Let $m \in \N$.  By the standard Morita equivalence $\Psi: {\rm M}_m(L(E)) \sim L(E)$ we have
the induced isomorphism $$\psi_m: K_0({\rm M}_m(L(E))) \rightarrow K_0(L(E))$$ for which
$\psi_m([1_{{\rm M}_m(L(E))}]) = m[1_{L(E)}]$.  By \cite[Proposition 9.3]{AT} there exists a
graph $M_mE$ for which $L(M_mE) \cong {\rm M}_m(L(E))$.  Specifically, this yields an
isomorphism $$\rho_m: K_0(L(M_mE)) \rightarrow K_0({\rm M}_m(L(E)))$$ for which
$\rho_m([1_{L(M_mE)}]) = [1_{{\rm M}_m(L(E))}]$.  Now the composition $$\kappa =
\rho_d^{-1} \circ \psi_d^{-1} \circ \varphi \circ \psi_c \circ \rho_c$$ is an isomorphism
from $K_0(L(M_cE))$  to $K_0(L(M_dE))$ for which $\kappa([1_{L(M_cE)}])=[1_{L(M_dE)}]$.

Since $L(M_cE)$ and $L(M_dE)$ are both Morita equivalent to $L(E)$ and therefore to each
other, the existence of the isomorphism $\kappa$ having the indicated property allows us to apply the
aforementioned algebraic Kirchberg Phillips Theorem \cite[Theorem 2.5]{ALPS}, from which we
conclude that there is an isomorphism $L(M_cE) \cong L(M_dE)$, which yields ${\rm M}_c(L(E))
\cong {\rm M}_d(L(E))$ as desired.

($\Leftarrow$) \   Conversely, suppose $\gcd(c,n) \neq \gcd(d,n)$.  In this case there
cannot be a ring isomorphism from ${\rm M}_c(L(E))$ to ${\rm M}_d(L(E))$, as otherwise,
by contradiction, if such exists then (by standard ring theory, see. e.g. \cite[p. 5]{Ros})
there would exist an isomorphism $$\tau: K_0({\rm M}_c(L(E))) \rightarrow K_0({\rm M}_d(L(E)))$$
for which $\tau([1_{{\rm M}_c(L(E))}]) = [1_{{\rm M}_d(L(E))}].$  But then
$\varphi = \psi_d \circ \tau \circ \psi^{-1}_c$ (with $\psi_m$ as above) would be an
automorphism of $G = K_0(L(E))$ for which $\varphi(c[1_{L(E)}]) = d[1_{L(E)}]$, which is
impossible by Lemma \ref{grouptheoryresult}.
\end{proof}

In fact, the proof given above yields that the converse direction of Theorem \ref{regularmodulefiniteorder} holds
for all rings $R$ for which $[1_R]$ has finite order in $K_0(R)$.

We note that Theorem \ref{regularmodulefiniteorder} generalizes \cite[Theorem 4.14]{AAP}
and \cite[Theorem 5.2]{AAP} (as well as \cite[Theorem 5.9]{AALP}) from the (purely infinite
simple) Leavitt algebras $L_q$ to all purely infinite simple unital Leavitt path algebras
for which $[1_{L(E)}]$ has finite order in $K_0(L(E))$, since
the order of $[1_{L_q}]$ in $K_0(L_q)$ is $q-1$.  In the related article \cite{AS2} we will
show that the indicated isomorphisms between matrix rings can be explicitly described.

To complete the determination of the Matrix Type of all purely infinite simple unital Leavitt
path algebras, we now consider the case where $[1_{L(E)}]$ has infinite order in $K_0(L(E))$.

\begin{lem}\label{eigenlemma}
Let $G$ be a finitely generated abelian group.  If there exists $m,n \in \N$, $\sigma \in {\rm Aut}(G)$, and
$x \in G$ of infinite order such that $n \sigma(x) = mx$, then $n = \pm m$.
\end{lem}
\begin{proof}
Since $G$ is finitely generated, $G \cong H \oplus \Z^t$ for some $t\in \N$, with $H$ finite.
By hypothesis, $x$ has nonzero component $\hat{x}$ in $\Z^t$.  We easily get that ${\rm Aut}(G) = {\rm Aut}(H) \oplus {\rm Aut}(\Z^t)$.  Since ${\rm Aut}(\Z^t) = GL(t,\Z)$,
since $ \sigma(\hat{x}) = \frac{m}{n}\hat{x}$ by hypothesis,  and since
the only rational eigenvalues of an invertible integer-valued matrix are $1$ and
$-1$, we have $\frac{m}{n} = \pm 1$, which gives the result.
\end{proof}

Following terminology introduced by P. Vamos, we say that a ring $R$ has {\it Invariant
Matrix Number} in case ${\rm M}_i(R) \not\cong {\rm M}_j(R)$ for every pair of positive
integers $i\neq j$.

\begin{prop}\label{IMTcorollary}
Let $R$ be a unital ring for which the order of $[1_{R}]$ in $K_0(R)$ is infinite, and for
which $K_0(R)$ is a finitely generated group.  Then $R$ has Invariant Matrix Number.

In particular, if $E$ is finite, and $[1_{L(E)}]$ has infinite order in
$K_0(L(E))$, then $L(E)$ has Invariant Matrix Number.
\end{prop}
\begin{proof}
Let $m,n \in \N$ and suppose ${\rm M}_m(R) \cong {\rm M}_n(R)$.  Then, as noted above, there
exists an isomorphism $\tau: K_0({\rm M}_m(R)) \rightarrow K_0({\rm M}_n(R))$ for which
$\tau([1_{{\rm M}_m(R)}]) = [1_{{\rm M}_n(R)}].$
 In addition, as noted previously, for any $a\in \N$, using
the standard Morita equivalence between $R$ and ${\rm M}_a(R)$  we get an isomorphism $\psi_a: K_0({\rm M}_a(R)) \rightarrow K_0(R)$  for which
$\psi_a([1_{{\rm M}_a(R)}])=a[1_{R}]$.
Then the composition
$$\sigma = \psi_n \circ \tau \circ \psi_m^{-1}: K_0(R) \rightarrow K_0(R)$$
is an automorphism of $K_0(R)$ for which  $m\sigma([1_{R}]) = n[1_{R}]$.  By Lemma
\ref{eigenlemma} and because $[1_R]$ has infinite order, $m = n$.

The result applies immediately to the indicated rings of the form $L(E)$ since,
as noted above, for these rings $K_0(L(E))$ is a finitely generated abelian group.
\end{proof}

The requirement that $K_0(R)$ be finitely generated cannot be removed from Proposition
\ref{IMTcorollary}.  We thank E. Pardo for providing the following
example and subsequent remarks.

\begin{example}
Let $R$ be any unital ring.    Consider the ring
$S=\varinjlim ({\rm M}_n(R),f_{m,n})$, where the connecting maps
$f_{m,n}:{\rm M}_m(R)\rightarrow {\rm M}_n(R)$ are defined when $m$ divides $n$, and
are the classical block diagonal maps.   It is well known (and not hard to show)
that $S\cong {\rm M}_n(S)$ for every
$n\in \mathbb{N}$.   (That is, $S$ has {\it Single Matrix Number}.)


Suppose also that $R$ has the property that $[1_R]$ has infinite order in $K_0(R)$.
 Since $K_0$ is a continuous functor, we have
$K_0(S)\cong \varinjlim K_0({\rm M}_n(R))$. As utilized above, we have
$(K_0({\rm M}_n(R)), [I_n])\cong (K_0(R), n[1_R])$. Since
$K_0(f_{m,n})([I_m])=[I_n]$, we conclude that the order of $[1_S]$ in
$K_0(S)$ is infinite as well.

Thus for $R$ any unital  ring for which $[1_R]$ has infinite order in $K_0(R)$,
the ring $S=\varinjlim ({\rm M}_n(R),f_{m,n})$ has the property that
 $[1_S]$ is of infinite order in $K_0(S)$, and for which $S$ does not have Invariant Matrix Number, as desired.
\end{example}

Our interest here is in purely infinite simple rings $R$, so one might ask whether the finitely generated hypothesis can be dropped from Proposition \ref{IMTcorollary} in case $R$ has this additional property.  But even in this case the finitely generated hypothesis on $K_0(R)$ is needed, since if one starts with $R$ purely infinite simple in the previous Example, then  $S$ can easily be shown to be purely infinite simple as well.

For any unital ring $R$, if we construct $S=\varinjlim {\rm M}_n(R)$ as in the
Example, then it is well known that
$$K_0(S) \cong \mathbb{Q}\otimes_\mathbb{Z} K_0(R).$$
(To establish this, for each $n\in \N$ and   $[A]\in K_0(M_n(R))$ define $\kappa_n:  K_0(M_n(R)) \rightarrow \mathbb{Q}\otimes_\mathbb{Z} K_0(R)$
by setting  $\kappa_n([A]) = \frac{1}{n} \otimes \psi_n([A])$.  One then verifies that the maps $\kappa_n$ are consistent with the maps $K_0(f_{m,n})$, which then yields a homomorphism from  $K_0(S)\cong \varinjlim K_0({\rm M}_n(R))$ to $\mathbb{Q}\otimes_\mathbb{Z} K_0(R).$  That this map is an isomorphism is easily shown by constructing the appropriate inverse map.)
In particular, for every $m,n \in \N$, we may define $\sigma \in Aut(K_0(S))$ as
the linear extension of $q \otimes t \mapsto \frac{m}{n}q \otimes t$ for $q \in \mathbb{Q}, \
 t\in K_0(R)$.   Then $\sigma \in {\rm Aut}(K_0(S))$, and for every
$x \in K_0(S)$ we have $n \sigma(x) = mx$.   (Compare this to the hypotheses of Lemma \ref{eigenlemma}.)
   This observation is
what accounts for the difference in the Matrix Type of the ring $S$ given here (i.e., Single Matrix Number), as compared to the Invariant Matrix Number
property of
 Leavitt path algebras of finite graphs $E$ for which $[1_{L(E)}]$ has infinite order in $K_0(L(E))$.

\singlespacing

\end{document}